\documentclass[12pt]{article}%

\makeatletter
\def\input@path{{styles/}}
\makeatother

\usepackage[bibencoding=utf8,style=alphabetic,backend=biber]{biblatex}%
\usepackage{sariel_biblatex}%

\usepackage[T1]{fontenc}

\usepackage{xspace}%
\usepackage[cm]{fullpage}%
\usepackage{amsmath}%
\usepackage{amssymb}%
\usepackage{mleftright}%

\usepackage{hyperref}%
\usepackage[ocgcolorlinks]{ocgx2}

\newcommand{\ProbC}{{\mathbb{P}}}
\newcommand{\ExC}{{\mathbb{E}}}

\newcommand{\Prob}[1]{\ProbC\mleft[ #1 \mright]}
\newcommand{\Ex}[1]{\ExC\mleft[ #1 \mright]}

\newcommand{\ExCond}[2]{\ExC\!\left[%
       #1 \;\middle\vert\; #2 \right]}

\hypersetup{%
      unicode,
      breaklinks,%
      colorlinks=true,%
      urlcolor=[rgb]{0.25,0.0,0.0},%
      linkcolor=[rgb]{0.5,0.0,0.0},%
      citecolor=[rgb]{0,0.2,0.445},%
      filecolor=[rgb]{0,0,0.4},
      anchorcolor=[rgb]={0.0,0.1,0.2}%
}

\newcommand{\myqedsymbol}{\rule{2mm}{2mm}}
\newcommand{\SarielThanks}[1]{%
   \thanks{%
      Department of Computer Science; %
      University of Illinois; %
      201 N. Goodwin Avenue; %
      Urbana, IL, 61801, USA; %
      \href{mailto:spam@illinois.edu}{sariel@illinois.edu}; %
      \url{http://sarielhp.org/}.%
   #1%
   }%
}

\usepackage[amsmath,thmmarks]{ntheorem}%
\theoremseparator{.}%

\theoremstyle{plain}%
\newtheorem{theorem}{Theorem}[section]

\theoremstyle{plain}%
\theoremheaderfont{\sf} \theorembodyfont{\upshape}%
\newtheorem*{remark:unnumbered}[theorem]{Remark}%
\newtheorem{remark}[theorem]{Remark}%

\theoremheaderfont{\em}%
\theorembodyfont{\upshape}%
\theoremstyle{nonumberplain}%
\theoremseparator{}%
\theoremsymbol{\myqedsymbol}%

\newcommand{\QuickSort}{\texttt{QuickSort}\xspace}

\bibliography{weak_chernoff}

\providecommand{\th}{th\xspace}%
\renewcommand{\th}{th\xspace}

\begin{document}

\title{An Easy Proof of a Weak Version of Chernoff Inequality}

\author{Sariel Har-Peled\SarielThanks{Work on this paper
      was partially supported by NSF AF award
      CCF-2317241.
   }}

\date{\today}

\maketitle

\begin{abstract}
    We prove an easy but very\footnote{More very{}s can be provided by
       the author upon request.}  weak version of Chernoff's
    inequality \cite{mr-ra-95, dp-cmara-09, m-18-fpcba,
       k-scapcbtg-25}. Namely, that the probability that in $6M$
    throws of a fair coin, one gets at most $M$ heads is $\leq 1/2^M$.
\end{abstract}

\section{The proof}

Consider a game that starts with $X_0 = 2^M$ balls. In the $i$\th
iteration, you flip a fair coin $C_i \in \{0,1\}$.  If $C_i=1$, which
happens with probability half, one throws away half the balls they
currently have. Otherwise, all balls are kept to the next iteration.
Let $Y_i$ be the number of balls you have after the $i$\th iteration,
and observe that
\begin{equation*}
    \ExCond{Y_i}{Y_{i-1}}
    =%
    Y_{i-1}/2 + Y_{i-1}/4 = (3/4)Y_{i-1}.
\end{equation*}
As such, we have
\begin{equation*}
    \Ex{Y_i}
    =%
    \Ex{\ExCond{Y_i}{Y_{i-1}} \bigr.}
    =%
    \Ex{(3/4)Y_{i-1}}
    =
    (3/4)^i 2^M.
\end{equation*}
Observe that $(3/4)^3 = 27/64 \leq 1/2$. Thus, we have
\begin{math}
    \Ex{Y_{6M}} = (3/4)^{6M} 2^M \leq 2^{M-2M } = 1/2^M.
\end{math}
By Markov's inequality, we have that
\begin{equation*}
    \Prob{ \sum\nolimits_{i=1}^{6M} C_i \leq M}
    =%
    \Prob{ Y_{6M} \geq 1 }
    \leq
    \frac{\Ex{Y_{6M}}}{1}
    \leq
    \frac{1}{2^M}.
    \hspace{3cm}\tag*{\myqedsymbol}
\end{equation*}

\begin{remark}
    The above proof is essentially Azuma's inequality proof without
    optimizing the constants. The magic is in using the formula
    $\Ex{\ExCond{X}{Y} \bigr.} =\Ex{X}$.
\end{remark}
\begin{remark}[Connection to \QuickSort]
    The proof came about from analyzing \QuickSort. Indeed consider an
    element $e$ in an array of size $n$ being sorted by \QuickSort,
    and let $Y_i$ be the size of the recursive subproblem containing
    $e$ in the $i$\th level of the recursion. A careful
    analysis\footnote{Here is the not careful analysis. With
       probability half the pivot for the subproblem of $Y_i$
       ``succeeds'' to fall in the range
       $[\tfrac{1}{4}Y_{i-1}, \tfrac{3}{4}Y_{i-1}]$. Then the size of
       the larger subproblem is $(3/4)Y_{i-1}$. Otherwise, we get
       nothing, and we pretend that $Y_i = Y_{i-1}$. Thus,
       $\ExCond{Y_i}{Y_{i-1}} \leq (7/8)Y_{i-1}$.}  shows that
    $\ExCond{Y_i}{Y_{i-1}} \leq (3/4)Y_{i-1}$. The same proof now goes
    through with $M=c\log_2 n$, to argue that the probability
    \QuickSort recursion depth exceeds $6M$ is at most
    $n \cdot 1 /2^{M} = 1/n^{c-1}$.
\end{remark}

\printbibliography

\end{document}